%% file: main.tex
\documentclass{article}
\usepackage{amsmath}[1996/11/01]
\usepackage{amssymb,amsthm,amsxtra}
\usepackage{epsfig}

\newtheorem{thm}{Theorem}

\newtheorem{rem}[thm]{\it{Remark}}

\numberwithin{equation}{section} \numberwithin{thm}{section}

\input{newcommands}

\begin{document}

\title{{\bf Riemann-Hilbert problems for last passage percolation}}
\author{{\bf Jinho Baik} \\ Princeton University and \\
Institute for Advanced Study}

\date{July 1, 2001}
\maketitle

\begin{abstract}
Last three years have seen new developments
in the theory of last passage percolation, which has
variety applications to
random permutations, random growth
and random vicious walks.
It turns out that a few class of models have determinant formulas for the
probability distribution, which can be analyzed asymptotically.
One of the tools for the
asymptotic analysis has been the Riemann-Hilbert method.
In this paper, we survey
the use of Riemann-Hilbert method in the last passage
percolation problems.
\end{abstract}

\section{introduction}\label{sec-intro}

Let $\Sigma$ be the unit circle in $\C$, oriented counterclockwise,
and let $\Omega_+=\{z\in\C : |z|<1\}$,
$\Omega_-=\{z\in\C : |z|>1\}$.
Set
\begin{equation}\label{e-varphi}
  \varphi(z)=e^{t(z+z^{-1})}.
\end{equation}
Consider the following Riemann-Hilbert problem :
$Y(z)$ is the $2\times 2$ matrix-valued function satisfying
\begin{equation}\label{e-rhp}
\begin{cases}
  Y(z) \quad \text{is analytic in $z\in\Omega_\pm$,
and is coninous in $\overline{\Omega_\pm}$},\\
  Y_+(z)= Y_-(z) \begin{pmatrix} 1& z^{-k} \varphi(z) \\
0&1 \end{pmatrix}, \quad z\in\Sigma, \\
  Y(z)= (I+O(z^{-1)}) z^{k\sigma_3}, \quad \text{as $z\to\infty$}.
\end{cases}
\end{equation}
Here $\sigma_3= \bigl( \begin{smallmatrix} 1&0\\0&-1\end{smallmatrix}
\bigr)$, and $Y_+(z)$ (resp., $Y_-(z)$) is the limit from the inside
(resp., outside) of $\Sigma$.

This Riemann-Hilbert problem appears in a class of problems which have many different
interpretations like longest increasing subsequence in a random permutation,  last passage
percolation, polynuclear growth model, and random vicious walks. In this paper, we survey a
list of problems related to the above Riemann-Hilbert problems and discuss various aspects
of the RHP \eqref{e-rhp}.

Before we state the problems, we first consider the solution to \eqref{e-rhp}.
From the work of Fokas, Its and Kitaev \cite{FIK}, the RHP \eqref{e-rhp}
is related to the orthogonal polynomials on the unit circle.
Indeed (see e.g, \cite{BDJ}),
if we let $\pi_k(z)$ be the $k^{th}$ monic (leading coefficient=1) orthogonal
polynomial with respect to the measure $\varphi(z)dz/(2\pi iz)$
on $\Sigma$,
\begin{equation}\label{e-opN}
  \int_{|z|=1} \pi_k(z)\overline{z^j} \varphi(z)\frac{dz}{2\pi iz} =
N_k\delta_{jk}
\qquad 0\le j\le k,
\end{equation}
the function
\begin{equation}\label{e-Y}
  Y(z)= \begin{pmatrix}
\pi_k(z) & (C\pi_k)(z) \\
-N_{k-1}^{-1}\pi^*_{k-1}(z) & -N_{k-1}^{-1} (C\pi^*_{k-1})(z)
\end{pmatrix}
\end{equation}
is a solution to \eqref{e-rhp},
where $\pi^*_{n}(z)=z^{n}\pi_{n}(1/z)$, and $(Cf)(z)$ is the Cauchy transform
of $f$ :
\begin{equation}
  (Cf)(z)=\frac1{2\pi i} \int_{\Sigma} \frac{f(s)}{s-z} ds.
\end{equation}
Also it is a standard argument to show that the solution is unique. Here $\varphi$ can be a
general function, not necessarily of the form \eqref{e-varphi}. Thus for instance,
\begin{equation}\label{e-N}
   \pi_k(z)= Y_{11}(z;k), \qquad N_{k-1}^{-1}= - Y_{21}(0;k).
\end{equation}

In this paper, we discuss three different uses of this RHP formulation for orthogonal
polynomials $\pi_k$ : (1) obtain asymptotic result using the Deift-Zhou steepest-descent
method (2) obtain differential (or difference) equations (3) make a connection to so-called
integrable operators. These three topics are discussed in Section \ref{sec-asym}, Section
\ref{sec-diff} and Section \ref{sec-intg}, respectively. The rest of the paper is organized
as follows. In Section \ref{sec-perc}, we state the last passage percolation problems whose
solution has Toeplitz/Hankel determinant formulas, which can be expressed in terms of the
solution to the RHP \eqref{e-rhp}. These percolation problems have various different
interpretations, and Section \ref{sec-appl} discusses some of these applications. There are
a few other last passage percolation problems which can also be solved in terms of the RHP
of the type \eqref{e-rhp}, but with a different function $\varphi$. These problem are
discussed in the final section \ref{sec-othe}. Most of the material in this paper have
appeared somewhere. We indicate the references in each section. The difference equation,
Theorem \ref{thm-diff} is the only new result.

\medskip
\noindent {\bf Acknowledgments.} The author would like to thank
Xin Zhou and Ken McLaughlin for inviting him to the
AMS special session on Riemann-Hilbert problems held in
Birmingham, AL.
The author would also like to thank Anne Boutet de Monvel for
kindly inviting him to Universit\'e Paris 7,
where a part of this paper is written.
This work was supported in part by NSF Grant \# DMS 97-29992.

\section{Last passage percolation}\label{sec-perc}

\subsection*{Square case}

\begin{figure}[ht]
 \centerline{\epsfig{file=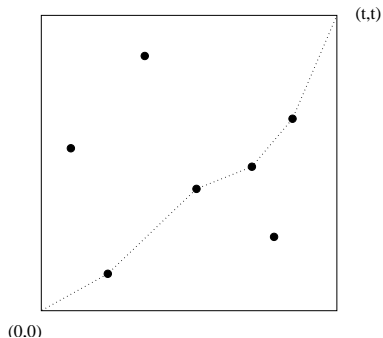, width=5cm}}
 \caption{Poisson points in $(0,0)\times (t,t)$ and the longest up/right path}
\label{fig-square}
\end{figure}

Consider a Poisson process of rate $1$ in the plane. An up/right path is, by definition, a
piece-wise linear curve with positive slope where defined. The `length' of an up/right path
is defined by the number of Poisson points on the path. Let $L(t)$ be the length of the
`longest' up/right path starting from $(0,0)$ ending at $(t,t)$ (see Figure
\ref{fig-square}).
See Section \ref{sec-appl} below for
other different interpretations.
The basic formula for us is the following,
proved by \cite{Gessel} (see also \cite{Rains}) :
\begin{equation}\label{e-U}
\Prob(L(t)\le \ell)=\frac1{Z} \Exp_{U\in U(\ell)}(\det(\varphi(U)))
\end{equation}
where $\varphi$ is given in \eqref{e-varphi}
and $Z=e^{t^2}$ is the normalization
constant to make the left hand side $1$ as $\ell\to\infty$.

As is well-known, using the Weyl's integration formula
for the unitary group, the above expected value is equal to
a Toeplitz determinant :
\begin{equation}\label{e-Toepl}
  \Exp_{U\in U(\ell)} \det(\varphi(U))
= \det(\varphi_{j-k})_{0\le j,k<\ell} =:
D_\ell(t),
\end{equation}
where $\varphi_j$ is the $j^{th}$ Fourier coefficient of
$\varphi$.
Now we will relate this Toeplitz determinant with the RHP \eqref{e-rhp}.

As in Section \ref{sec-intro}, let $\pi_k(z)$ be the
$k^{th}$ monic orthogonal polynomial with respect to $\varphi(z)dz/(2\pi iz)$
on the unit circle and $N_k$ be defined by \eqref{e-opN}.
It is direct to check that the orthogonal polynomials
have the determinant expression (see e.g., \cite{Szego})
\begin{equation}\label{e-opdet}
  \pi_k(z):= \frac1{D_k(\varphi)} \det\begin{pmatrix}
\varphi_0 & \varphi_{1} & \cdots & \varphi_{k} \\
\varphi_{-1} & \varphi_0 & \cdots & \varphi_{k-1} \\
\vdots & \vdots & \vdots &\vdots \\
\varphi_{-k+1} & \varphi_{-k+2} & \cdots & \varphi_{1}\\
1 & z & \cdots & z^k
\end{pmatrix}.
\end{equation}
Thus we have
\begin{equation}\label{e-ND}
  N_k= \frac{D_{k+1}}{D_k}.
\end{equation}
Note from \eqref{e-U} and \eqref{e-Toepl} that
\begin{equation}
  e^{t^2}=Z= \lim_{n\to\infty} D_{n},
\end{equation}
which also can be seen from the strong Szeg\"o limit theorem.
Hence we have from \eqref{e-N}
\begin{equation}\label{e-ratioD}
  \Prob(L(t)\le \ell) = e^{-t^2}D_\ell=
\lim_{n\to\infty} \frac{D_\ell}{D_n}
= \prod_{k=\ell}^\infty N_k^{-1}
= \prod_{k=\ell}^\infty ( -Y_{21}(0; k+1)) .
\end{equation}

\begin{rem}\label{rem-daku}
  In addition to the Toeplitz determinant formula, there are two other
determinant formulas for $\Prob(L(t)\le \ell)$.
Both of them are Fredholm determinants of operators, one acting on
$\Sigma$ (see \cite{BDJ2}),
and the other acting on the discrete set $\{\ell, \ell+1, \cdots\}$
(see \cite{BOO, kurtj:disc, Ok99}).
The first Fredholm determinant
has the RHP expression which is algebraically equivalent
to \eqref{e-rhp}.
See Section \ref{sec-intg} below for details.
The second Fredholm determinant also has the RHP expression,
but now the RHP has a jump
condition on the discrete set.
However, for the second formula,
the kernel of the operator has integral representation, and hence
the classical steepest-descent method is sufficient to obtain
the asymptotics (\cite{BOO, kurtj:disc}).
Thus the limit of $\Prob(L(t)\le \ell)$ can be obtained
from the second algebraic formula, which avoids the RHP analysis
in Section \ref{sec-asym} below.
But for some applications
(see \cite{kurtj:trans, seppal})
one also needs uniform tail estimates of the `scaled' random variable
$(L(t)-2t)/t^{1/3}$
in $t$ (see \eqref{e-est1}, \eqref{e-est2} below).
This problem is more difficult than the convergence in distribution,
and the lower bound \eqref{e-est2} has been obtained only from RHP analysis
so far \cite{BDJ}.
The identities between the Toeplitz determinants and the Fredholm
determinants are discussed in \cite{BDJ2, BoOk, BasorW, BDR, Bott1, Bott2}.
\end{rem}

\subsection*{Triangle case}

\begin{figure}[ht]
 \centerline{\epsfig{file=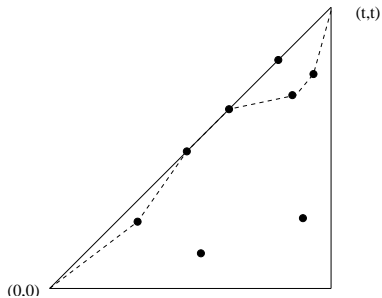, width=5cm}}
 \caption{Triangle case}
\label{fig-triangle}
\end{figure}

Suppose now that we take a (2-dimensional)
Poisson process of rate 1 in the half plane $y<x$,
and take a (1-dimensional)
Poisson process of rate $\alpha \ge 0$ on the line $y=x$
(see Figure \ref{fig-triangle}).
Let $L_s(t;\alpha)$ be again the length of the longest up/right path
from $(0,0)$ to $(t,t)$.
The algebraic formula for this case is obtained in \cite{Rains}
for $\alpha=0$ and in \cite{BR1} for general $\alpha>0$ :
\begin{equation}\label{e-O}
  \Prob(L_s(t;\alpha) \le \ell)
= \frac1{Z} \Exp_{U\in O(\ell)} \det( (1+\alpha U) e^{tU}),
\end{equation}
where $Z=e^{\alpha t+t^2/2}$.

There are two components of $O(\ell)$ depending on the sign
of the determinant of the matrix.
Using the Weyl's integration formula,
the expected value over each component is expressed in terms of one of
the Hankel determinants of the form (see e.g., Theorem 2.2 of \cite{BR1})
$\tilde{H}_\ell=\det(h_{j+k})_{\ell\times \ell}$
where
\begin{equation}
  h_j=\int_{-1}^1 x^j h(x) (1-x)^{\pm\frac12}(1+x)^{\pm\frac12} dx,
\end{equation}
with
\begin{equation}\label{e-h}
 h(x)= (1+\alpha^2+2\alpha x) e^{2tx}
= (1+\alpha z)(1+\alpha /z) \varphi(z),
\qquad x=\frac12(z+z^{-1}).
\end{equation}
As an analogue of \eqref{e-opdet},
these Hankel determinants are related to the orthogonal polynomials
on the interval $(-1,1)$.
But the orthogonal polynomials on the unit circle
with weight $(1+\alpha z)(1+\alpha /z) \varphi(z)$
and the orthogonal polynomials on the interval $(-1,1)$
with weight $h(x)(1-x)^{\pm\frac12}(1+x)^{\pm\frac12}$
are related in a simple way (see \cite{Szego}).
Thus the ration of Hankel determinants $\tilde{H}_{\ell}/\tilde{H}_{\ell+1}$
can be expressed (see theorem 2.3 of \cite{BR1})
in terms of orthogonal polynomials
(norms and the values at $z=0$) on the unit circle
$\Sigma$ with respect to the measure
\begin{equation}\label{e-tem}
 (1+\alpha z)(1+\alpha /z) \varphi(z) \frac{dz}{2\pi iz}.
\end{equation}
Moreover, after some algebraic work (see Section 3 of \cite{BR1}),
the simple factor
$(1+\alpha z)(1+\alpha /z)$ in the above measure
can be removed and $\tilde{H}_{\ell}/\tilde{H}_{\ell+1}$
can be expressed purely in terms of orthogonal polynomials $\pi_k$
for the measure $\varphi(z)dz/(2\pi iz)$.
The result is \cite{BR1}
\begin{equation}
  \Prob(L_s(t;\alpha)\le 2\ell+1)
= e^{-\alpha t} \frac12\bigl\{
(\pi^*_{2\ell}(-\alpha) + \alpha \pi_{2\ell}(-\alpha))
H_{\ell}^+ + (\pi^*_{2\ell}(-\alpha) - \alpha \pi_{2\ell}(-\alpha))
H_{\ell}^- \bigr\}
\end{equation}
where
\begin{equation}\label{e-Haneklell}
  H^\pm_{\ell} = \prod_{k=\ell}^\infty N_{2k+1}^{-1} ( 1\mp \pi_{2k+1}(0)).
\end{equation}
There is a similar formula for $\Prob(L_t(t;\alpha)\le 2\ell)$.
Thus again, the probability distribution for $L_s(t;\alpha)$ can
be expressed in terms of the solution of the RHP \eqref{e-rhp}.

\begin{rem}\label{rem-symm}
  In \cite{BR1, BR2, BR3}, the authors considered five types of symmetry
of the Poisson model.
Reflection symmetry about the anti-diagonal,
reflection symmetry about both diagonal and anti-diagonal,
and rotation symmetry about the center are considered
in addition to the square case and the triangle case.
The distribution of the
longest up/right path in each case has the determinantal expression
(either Toeplitz or Hankel determinants with a simple change of the weight),
and can be expressed in terms of the RHP \eqref{e-rhp}.
\end{rem}

\begin{rem}
 As in square case (see Remark \ref{rem-daku}),
there is also a different algebraic formula for the triangle case
\cite{Rains:corr}.
But the result is so-called Fredholm Pfaffian of an operator acting on
a discrete set.
The kernel of the operator is rather involved, and so far there has been
no work for the asymptotic analysis from this Fredholm Pfaffian formula.
\end{rem}

\subsection*{External sources}

\begin{figure}[ht]
 \centerline{\epsfig{file=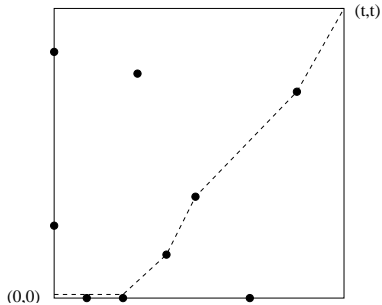, width=5cm}}
 \caption{External sources case}
\label{fig-external}
\end{figure}

Suppose that in the square model above, we have
additional (1-dimensional) Poisson processes of rate $\alpha_+$ and $\alpha_-$
on the lines $y=0$ and $x=0$, respectively.
We assume that the corner $(0,0)$ has no point
(see Figure \ref{fig-external}).
Let $L_e(t;\alpha_+,\alpha_-)$ be the length of the longest
up/right path from $(0,0)$ to $(t,t)$.
This problem arises from a polynuclear growth model with random initial data
\cite{SpohnP2}, \cite{SpohnP3}.

The distribution in this case is given in \cite{BR4} :
\begin{equation}\label{e-extform}
  \Prob(L_e(t;\alpha_+,\alpha_-)
= \frac1{Z} (D'_{\ell}-\alpha_+\alpha_-D'_{\ell-1})
\end{equation}
where
\begin{equation}
 D'_{\ell}:=
 \Exp_{U\in U(\ell)}
\det((1+\alpha_+U)(1+\alpha_- U^{-1}) \varphi(U)),
\end{equation}
and $Z=e^{\alpha_+ t+\alpha_- t+t^2}$.

By a similar argument as in \eqref{e-tem}, $D'_{\ell}$ can be
written as (see Theorem 3.2 of \cite{BR1})
\begin{equation}
  D'_\ell= \frac{\pi^*_{\ell}(-\alpha_+)\pi^*(-\alpha_-)
-\alpha_+\alpha_-\pi_\ell(-\alpha_+)\pi_\ell(-\alpha_-)}{1-\alpha_+\alpha_-}
D_\ell, \qquad \alpha_+\alpha_-\neq 1,
\end{equation}
where $D_\ell=D_\ell(t)$ is given in \eqref{e-Toepl}.
For $\alpha_+\alpha_-=1$, we use the l'Hopital's rule in the above formula.
Thus again the distribution for $L_e(t;\alpha_+,\alpha_-)$
can be expressed in terms of the solution $Y$ of the RHP \eqref{e-rhp}.

\begin{rem}
 In the triangular case above, we may put an additional (1-dimensional)
Poisson process of rate $\alpha_+$ on the line $y=0$.
Then there still is a determinant formula for the distribution
of the longest up/right
path, which again can be expressed in terms of the RHP \eqref{e-rhp}.
\end{rem}

\section{Applications}\label{sec-appl}

In this section, we discuss various different interpretations of the
last passage percolation problems in Section \ref{sec-perc}.

\subsection*{Random permutations}

Given a permutation $\pi$,
a subsequence
$\{\pi(i_1),\cdots, \pi(i_k)\}$ such that
$i_1<\cdots<i_k$ and $\pi(i_1)<\cdots \pi(i_k)$ is called `increasing'.
The length of such an increasing subsequence is defined to be $k$.
Let $L_N(\pi)$ be the length of the longest increasing subsequence
of a permutation $\pi\in S_N$.
Now take $\pi$ randomly from $S_N$, and we ask the distribution of $L_N$.
This problem is known as Ulam's problem since early 1960's
(see e.g., \cite{AD}, \cite{BDJ}).

If we take $N$ itself as a Poisson random variable of rate $t^2$,
the distribution of $L_N$ is same as the distribution of $L(t)$ above :
\begin{equation}
  \Prob(L(t)\le \ell) = \sum_{N=0}^\infty
\frac{e^{-t^2}(t^2)^N}{N!} \Prob(L_N\le \ell).
\end{equation}
In other words, $\Prob(L(t)\le \ell)$ is the exponential generating function
of $\Prob(L_N\le\ell)$.
This follows from the fact that a configuration of $N$ points
in a square can
be regarded as a permutation $\pi\in S_N$ : consider the relative $x$-orders,
then relative $y$-orders.
The longest up/right path in a points configuration is precisely
the longest increasing subsequence of the corresponding
permutation.
Thus the asymptotics of $L_N$ as $N\to\infty$ can be obtained
from the asymptotics of $L(t)$ as $t\to\infty$ (see \cite{Jo2}).

The limiting distribution of $L(t)$ (see \eqref{e-FGUE} below), and thus the
limiting distribution of $L_N$, is first obtained in \cite{BDJ} by
using the RHP \eqref{e-rhp}.
Later the same result was obtained by \cite{Ok, BOO, kurtj:disc} using
different methods
which avoids the RHP analysis.
The uniform tail estimates (see \eqref{e-est1}, \eqref{e-est2}) and the
convergence of moments of the scaled random variable were
obtained in \cite{BDJ}.
As mentioned in Remark \ref{rem-daku},
the (lower) tail estimate obtained in \cite{BDJ} was later used in
other applications \cite{kurtj:trans}, \cite{seppal}.

Now we consider the corresponding permutation model for the triangle case.
Given random points in the triangle model, take mirror image of the
points in $y<x$ about the line $y=x$.
Hence we have points in the whole plane which is symmetric about $y=x$.
It is each to see that in terms of permutation, the symmetry condition
implies that $\pi^2=1$.
Also the points on the line $y=x$ corresponds to the fixed points $\pi(i)=i$.
Thus it is natural to consider the set
$\tilde{S}_{m,n}=\{ \pi\in S_{2m+n} : \pi^2=1, |\{ x: \pi(x)=x\}|=n \}$.
We equip this set with the uniform probability measure and consider the
longest increasing subsequence of a random $\pi\in \tilde{S}_{m,n}$.
As before, the asymptotics of $L_s(t;\alpha)$ implies the asymptotics
of the longest increasing subsequence in this set. See \cite{BR1, BR2, BR3}
for details.

\subsection*{Plancherel measure}

Let $Y_N$ be the set of partitions (Young diagrams)
$\lambda=(\lambda_1\ge \lambda_2\ge \cdots\ge \lambda_N)$,
$\lambda_j\in\N\cup\{0\}$, $\sum_j \lambda_j =N$.
Given a partition $\lambda$, let $d_\lambda$ be the number of
standard Young tableaux of shape $\lambda$
(see e.g. \cite{Stl} for the definitions).
It is a basic fact in the representation of symmetric group that
$\sum_{\lambda\vdash N} d_\lambda^2=|S_N|=N!$.
Define the Plancherel measure on $Y_N$ by
\begin{equation}
  \Prob( \lambda)= \frac{d_\lambda^2}{N!}, \qquad \lambda\in Y_N.
\end{equation}
The famous Robinson-Schensted (RS) correspondence
(see e.g. \cite{Sc}, \cite{Stl}) implies that
the first part $\lambda_1$ from the
Plancherel measure has the same distribution as the longest
increasing subsequence of random permutation $\pi\in S_N$ from the uniform
measure.
Thus from the result of the longest increasing subsequence,
the asymptotics of the first row $\lambda_1$ can be obtained.
The convergence in distribution and the convergence of moments
for the second row $\lambda_2$ was obtained in \cite{BDJ2}.
For the general rows $\lambda_j$, $j\ge 3$, the convergence in distribution
is obtained in \cite{Ok, BOO, kurtj:disc}, and
the convergence of moments are obtained recently in \cite{BDR}.

Let $\beta\ge 0$.
As a generalization, consider the $\beta$-Plancherel measure on $Y_N$
defined by
\begin{equation}
  \Prob( \lambda)= \frac{d_\lambda^\beta}{\sum_{\mu\in Y_N} d_\mu^\beta},
\qquad \lambda\in Y_N.
\end{equation}
One can regards $\beta$ as the inverse temperature.
The above Plancherel measure corresponds to the case when $\beta=2$.
The $1$-Plancherel measure also has a simple permutation interpretation.
Consider the set of involutions $\tilde{S}_N:=\{\pi\in S_N : \pi^2=1\}$.
Equip $\tilde{S}_N$ with the uniform probability measure and consider
the length of the longest increasing subsequence
of a random $\pi\in \tilde{S}_N$.
The RS correspondence implies that the first row $\lambda_1$ of
$\lambda\in Y_N$ under the $1$-Plancherel measure has
the same distribution as the longest increasing subsequence of a
random involution.
In \cite{BR2} (see also \cite{BR3}), the limiting distribution
of the first row of the $1$-Plancherel measure was obtained.
It would be interesting to obtain the limiting distributions for the
general $\beta$.

\subsection*{Polynuclear growth model}

Consider a flat 1-dimensional substrate.
There are random nucleation events happening on the substrate,
which can be thought as a 2-dimensional Poisson process
in space-time coordinates.
An island of height 1 with width 0 is made at each nucleation position,
and the island grows laterally with speed 1.
When two islands meet, they just form one island and keep growing at
the edges.
Further nucleation events may happen on top of islands.
This model is called Polynuclear growth (PNG) model.
An interest is the height $h(x,t)$ as time $t$ tends to infinity.

It is an observation by Pr\"ahofer and Spohn that
for certain choices of initial data, PNG models are in bijection to
the last passage percolation models in Section \ref{sec-perc}.
As the first case, suppose that initially there is an island of width 0 at
the position $x=0$.
It grows with speed 1, and we assume that further nucleation events
occur only on top of the base island
and on top of islands on the base island.
Then the height $h(0,t)$ is equal, in the sense of distribution,
to $L(t/\sqrt{2})$ of the square case
(see \cite{SpohnP1}).
Secondly, we consider the half infinite case : $x\in [0,\infty)$.
When an island hits the edge $x=0$, it stops growing at that side.
As before there is the base island at the position $x=0$
at time $t=0$.
In addition to the nucleation events of rate 1 in $(0,\infty)$,
and we assume that there are extra nucleation events at the position $x=0$
of rate $\alpha$.
Then the height $h(0,t)$ is equal to $L(t/\sqrt{2};\alpha)$ of
the triangle case (see \cite{SpohnP3}).
If one considers flat case, i.e. no base island, the PNG model is
in bijection to yet a different symmetry model of the Poisson process
mentioned in Remark \ref{rem-symm} above (see \cite{SpohnP1}).
Finally the case of external sources in Section \ref{sec-perc}
has the PNG model which has random initial data.
See \cite{SpohnP2, SpohnP3} for more details.

\subsection*{Random turn vicious walks}

At time $t=0$, each site of the semi-infinite lattice $\{\cdots, -2,-1,0\}$
is occupied by a walker.
A walker is called `right-movable' if its
right-neighboring site is vacant.
At each discrete time, we pick a random walker among
all right-movable walkers (there are finitely many such walkers), and
move it to its right-neighboring site.
Hence at each time, one and only one walker moves to its right.
We let the walkers walk by the above rule for $N$ time steps.
In the next $N$ time steps, we let walkers move to their \emph{left}
according to the similar rule.
We further impose the condition that after total $2N$ time steps,
all the walkers return to their original positions.
The interest is the position $x_1(N)$
of the right-most walker at time $N$.
It is shown in \cite{forrester99} (see also \cite{GOV}) that
the distribution of $x_1(N)$ is same as the distribution of $L_N$, the
length of the longest increasing subsequence of a random permutation.

If we just let the walkers walk up to time $N$ without assuming that
they all return to their original positions, $x_1(N)$ is,
in the sense of distribution,
equal to the length of the longest increasing subsequence of a random
involution ($\pi^2=1$).
See \cite{BR3} for the asymptotic results of this case.

\subsection*{Queueing theory and Totally asymmetric simple exclusion process}

Lattice versions of the last passage percolation problems are known to be
related to the queueing theory and also so-called totally asymmetric simple
exclusion processes (TASEP) (see e.g. \cite{kurtj:shape}).
See Section \ref{sec-othe} below for lattice last passage percolation
problems.

\section{Asymptotics}\label{sec-asym}

In this section, we sketch asymptotic analysis of the RHP \eqref{e-rhp}
and discuss the asymptotic results for the random variables
$L(t)$, $L_s(t;\alpha)$ and $L_e(t;\alpha_+,\alpha_-)$.

The asymptotics of the RHP \eqref{e-rhp}
as $k,t\to\infty$ in all different regimes
was first considered
in \cite{BDJ}, and later in a little more detail in \cite{BR2}.
A special interest is the case when $k$ and $t$ are related as
\begin{equation}\label{e-scaling}
  \frac{2t}{k}= 1-\frac{x}{2^{1/3}k^{2/3}}.
\end{equation}
For $x$ in a compact set, the RHP \eqref{e-rhp}
is localized around $z=-1$ as $k\to\infty$,
and the local RHP is, after scaling,
the RHP for the Painlev\'e II solution.
In the below, we will discuss, heuristically, how
the Painlev\'e II solution arises.
This heuristic argument is taken from Section 6 of \cite{BDJ2}.

First, we show why the point $z=-1$ plays an important role
in the asymptotic analysis.
Since the RHP \eqref{e-rhp} is not normalized at $\infty$, we
algebraically transform it to a normalized problem.
Set
\begin{equation}\label{e-mdef1}
  m(z):= \begin{pmatrix} 0&-1\\1&0 \end{pmatrix} Y(z)
\begin{pmatrix} e^{tz}&0\\0&e^{-tz} \end{pmatrix}  , \qquad |z|<1,
\end{equation}
and set
\begin{equation}\label{e-mdef2}
  m(z):= \begin{pmatrix} 0&-1\\1&0 \end{pmatrix} Y(z)
\begin{pmatrix} 0&z^{-k}e^{tz^{-1}}\\-z^ke^{-tz^{-1}}&0\end{pmatrix},
\qquad |z|>1.
\end{equation}
Then $m(z)$ solves the new Riemann-Hilbert problem
\begin{equation}\label{e-newrhp}
\begin{cases}
  m(z) \qquad \text{is analytic in $z\in \Omega_\pm$,
and is continous in $\overline{\Omega_\pm}$,}\\
  m_+(z)=m_-(z) V(z), \qquad z\in\Sigma,\\
  m(z) = I +O(z^{-1}), \qquad z\to\infty,
\end{cases}
\end{equation}
where
\begin{equation}
V(z)= \begin{pmatrix} 0&-z^{-k}\psi(z)^{-1}\\
z^{k}\psi(z)&1  \end{pmatrix},
\qquad \psi(z):= e^{t(z-z^{-1})}.
\end{equation}
Now the oscillatory factors in the jump matrix are $V_{12}$ and
$V_{21}=-1/V_{12}$.
It is direct to check that the critical points of
\begin{equation}
  \log V_{21}(z) = k ( \frac{t}{k}(z-z^{-1})+\log z)
\end{equation}
are
\begin{equation}
  z= -\frac{k}{2t} \pm \sqrt{(\frac{k}{2t})^2-1}.
\end{equation}
When $2t<k$, there are two real critical points.
These two points collapses to one point $z=-1$ when $2t=k$.
And when $2t>k$, there are two complex critical points.
Thus one can imagine that the nature of asymptotics changes at $2t=k$.
This is indeed the case.
It is well-known that in the (so-called Plancherel-Rotach type)
asymptotics of the orthogonal polynomials,
the equilibrium measure of an associated variational problem
plays a key role (see e.g. \cite{DKMVZ2, DKMVZ3}).
For the case at hand, the support of the equilibrium measure
is the full circle when $2t\le k$,
and it is a part of the circle when $2t>k$.
And the point $z=-1$ is the place where the `gap' of the support
starts to open up at $2t=k$.

Thus it is natural to analyze the problem near $z=-1$
and in the regime $2t\sim k$.
The scaling \eqref{e-scaling} is chosen so that
the local RHP around $z=-1$ becomes non-trivial.
Writing $z=-1+s$ for $z$ near $-1$, we have
\begin{equation}
  \log V_{12}(z) -\pi i= -\frac{x}{k^{1/3}}(k^{1/3}s)
- \frac{x}{2^{4/3}k^{1/3}} (k^{1/3}s)^2
+ \frac{(k^{1/3}s)^3}{6} (1- \frac{3x}{2^{1/3}k^{2/3}} )
+ \cdots.
\end{equation}
If we take $w= 2^{-4/3}k^{1/3}s$, then we have
\begin{equation}
  \log V_{12}(z) -\pi i= -2xw + \frac83 w^3 + \cdots,
\end{equation}
and hence the jump matrix near $z=-1$ becomes as $k\to\infty$,
\begin{equation}
  \tilde{V}(w) = \begin{pmatrix} 0&-(-1)^k e^{-2(-xw+\frac43w^3)} \\
(-1)^k e^{2(-xw+\frac43w^3)}&1 \end{pmatrix}
\end{equation}
for $w\in i\R$.
But after rotation by $\pi/2$ and
removing $(-1)^k$ by a simple conjugation,
this is precisely the RHP for Painlev\'e II
equation with the choice of parameters $p=-q=1$, $r=0$ (see e.g. \cite{DZ2},
\cite{BDJ}).
Of course, this is only a heuristics, and one needs to justify
the convergence of the original RHP \eqref{e-rhp} to the
above local RHP.
As mentioned above, the justification proceeds differently
for $2t<k$ and $2t>k$, due to the change of the support of the equilibrium
measure. See \cite{BDJ}, \cite{BR2} for details.

The asymptotic results we obtain are, for example, the following
\cite{BDJ, BR2} :
There are numerical constants $C,c, x_0>0$ such that for large $k$ and $t$,
\begin{enumerate}
 \item if $x\ge x_0$,
\begin{equation}\label{e-asy1}
  \biggl|-Y_{21}(0;k)-1\biggr|, \qquad
\biggl|Y_{11}(0;k)\biggr| \le \frac{C}{k^{1/3}}e^{-cx^{3/2}},
\end{equation}
\item if $-x_0\le x\le x_0$,
\begin{equation}\label{e-asy2}
  \biggl|-Y_{21}(0;k)-1-\frac{2^{1/3}}{k^{1/3}}v(x)\biggr|,  \qquad
\biggl|Y_{11}(0;k)+(-1)^k\frac{2^{1/3}}{k^{1/3}}u(x)\biggr|
\le \frac{C}{k^{2/3}},
\end{equation}
 \item if $x\le -x_0$,
\begin{equation}\label{e-asy3}
  \biggl|-\sqrt{\frac{2t}{k}}e^{k(2t/k-\log (2t/k)-1)}Y_{21}(0;k)-1\biggr|,
\qquad
\biggl|(-1)^k\frac{2t}{2t-k}Y_{11}(0;k)\biggr| \le \frac{C}{2t-k}.
\end{equation}
\end{enumerate}
Here $v(x)$, $u(x)$ in case (ii) are given by
the Painlev\'e II solution
\begin{equation}
\begin{cases}
  u''=2u^3+xu, \\
u(x)\sim -\frac1{2\sqrt\pi x^{1/4}}e^{-2/3x^{3/2}}, \quad x\to +\infty.
\end{cases}
\end{equation}
and
\begin{equation}
 v(x)=\int_\infty^x (u(s))^2 ds.
\end{equation}
We also obtain the asymptotics of $Y(z)$ for general points $z\in\C$
\cite{BR2}.

Now using \eqref{e-asy1}-\eqref{e-asy3},
\eqref{e-ratioD} and \eqref{e-N} yield the convergence in distribution
in the square case \cite{BDJ} :
for fixed $x\in\R$,
\begin{equation}\label{e-FGUE}
 \lim_{t\to\infty}  \Prob\biggl( \frac{L(t)-2t}{t^{1/3}} \le x \biggr) =
\exp \biggl( -\int_x^\infty (y-x)(u(y))^2 dy \biggr)
=: F_{GUE}(x).
\end{equation}
The function $F_{GUE}(x)$ in \eqref{e-FGUE}
is called the GUE Tracy-Widom distribution \cite{TW1},
which was first obtained as the
limiting distribution for the fluctuation of the largest eigenvalue
of a random matrix from the Gaussian unitary ensemble in the random matrix
theory (see e.g. \cite{Mehta} for an introduction to random matrix theory).
Also the uniform tail estimates can be obtained \cite{BDJ} :
there are constants $C, c, x_0>0$ such that for all $t>0$,
\begin{equation}\label{e-est1}
  \Prob\biggl( \frac{L(t)-2t}{t^{1/3}} \le x \biggr) \le Ce^{-cx^{3/2}},
\qquad x\ge x_0,
\end{equation}
and
\begin{equation}\label{e-est2}
  \Prob\biggl( \frac{L(t)-2t}{t^{1/3}} \le x \biggr) \le Ce^{-cx^{3}},
\qquad x\le -x_0.
\end{equation}
The upper estimate \eqref{e-est1} also follows from a
large deviation result of \cite{Se98}.
But the lower estimate \eqref{e-est2} is obtained only from the
RHP analysis so far.
These estimates imply the convergence of moments of the scaled
random variable $(L(t)-2t)/t^{1/3}$ \cite{BDJ}.
These estimates have also been crucial for the analysis of
the transversal fluctuation of the
longest up/right path \cite{kurtj:trans} and of the
perturbation of the equilibrium measure for a related dynamical system,
called stick process \cite{seppal}.

The asymptotic analysis of the RHP \eqref{e-rhp} also yields similar results
for the triangle case $L_s(t;\alpha)$ \cite{BR2}
and the external sources case $L_e(t;\alpha_+,\alpha_-)$ \cite{BR4}.
The results show more interesting feature.
Depending on the values $\alpha, \alpha_+,\alpha_-$,
the scaled random variables converge to different distributions.
For example, $(L_s(t;\alpha)-2t)/t^{1/3}$
converges in distribution to the so-called GSE TW-function \cite{TW2},
\begin{equation}
   F_{GSE}(x):=\frac12 \biggl(
e^{\frac12\int_x^\infty u(y)dy }
+  e^{ -\frac12\int_x^\infty u(y)dy} \biggl) (F_{GUE}(x))^{1/2}
\end{equation}
for fixed $0\le \alpha<1$, and to the GOE TW-function \cite{TW2},
\begin{equation}
   F_{GOE}(x):= e^{ \frac12\int_x^\infty u(y)dy }
(F_{GUE}(x))^{1/2}
\end{equation}
for $\alpha=1$.
When $\alpha>1$ and fixed, a differently scaled random variable converges
to the Gaussian distribution.
Moreover if we scale $\alpha=1-\frac{2w}{t^{1/3}}$ with fixed $w\in\R$
and take $t\to\infty$, the above scaled random variable
converges to yet another distribution function for each fixed $w\in\R$.
This new one-parameter family of distributions interpolates
$F_{GSE}$ and $F_{GOE}$ as $w=\infty$ and $w=0$.
Similar feature appears for the analysis of
$L_e(t;\alpha_+,\alpha_-)$ where two-parameter family of new distributions
are obtained \cite{BR4}.

\section{Difference equation}\label{sec-diff}

In this section, we show that the RHP \eqref{e-rhp} yields difference
equations for the entries of $Y(0)$.
It is direct to check from \eqref{e-Y} that $Y_{11}(0)=-Y_{22}(0)$.
Also since the determinant of the jump matrix for \eqref{e-rhp}
is $1$, we have $\det Y(0)=1$.

\begin{thm}\label{thm-diff}
  Let
\begin{equation}
  Y(0;k)= \begin{pmatrix} -b(k)& d(k) \\ a(k)&b(k) \end{pmatrix}.
\end{equation}
Then $b$ satisfy the discrete Painlev\'e II equation
\begin{equation}\label{e-diffeq}
   \frac{k}{t}b(k)+ \bigl( b(k-1)+b(k+1) \bigr) (1-b(k)^2)=0 ,
\end{equation}
and $a$ and $d$ satisfy
\begin{eqnarray}
\label{e-diffeqa}   a(k)&=& \bigl(1-b(k)^2\bigr) a(k+1) , \\
\label{e-diffeqd}   d(k)&=& \bigl(1-b(k)^2\bigr) d(k-1).
\end{eqnarray}
\end{thm}

\begin{proof}
Consider \eqref{e-newrhp}.
Note that the jump matrix has determinant $1$.
Thus $\det(m_+)=\det(m_-)$, and hence $\det(m)$ is an entire function.
Also $\det(m)\to 1$ as $z\to\infty$.
Thus by Liouville's theorem, we have $\det(m(z))\equiv 1$, and in particular,
$(m(z))^{-1}$ is analytic in $\Omega_\pm$, and is continuous in
$\overline{\Omega_\pm}$.

Set
\begin{equation}\label{e-mas1}
   m(z)=m(z;k)= I+ \frac{m^\infty_1(k)}{z} + O(z^{-2}),
\qquad z\to\infty,
\end{equation}
and set $m(0;k)=B(k)$ and
\begin{equation}\label{e-mas2}
   m(z)=m(z;k)= B(k)(I+ m^0_1(k)z + O(z^{2})),
\qquad z\to 0.
\end{equation}

Set
\begin{equation}\label{e-Psi}
\Psi(z)=\Psi(z;k):= m(z)  \begin{pmatrix} 1&0\\0&z^k\psi(z) \end{pmatrix}  ,
\qquad
J_1:=   \begin{pmatrix} 0&0\\0&1 \end{pmatrix}  ,
\quad
J_0:=   \begin{pmatrix} 1&0\\0&0 \end{pmatrix}
\end{equation}
Note that by \eqref{e-mas1} and \eqref{e-mas2},
\begin{eqnarray}\label{e-Psias1}
  \Psi(z;k) &=& \bigl(I + \frac{m^\infty_1(k)}{z} + \frac{m^\infty_2(k)}{z^2}
+  O(z^{-3})\bigr)
e^{tzJ_1+k\log z J_1 -\frac{t}{z}J_1},
\quad z\to\infty, \\
 \label{e-Psias2}
\Psi(z;k) &=& B(k)\bigl(I + m^0_1(k)z + O(z^{2})\bigr)
e^{-\frac{t}{z}J_1+k\log z J_1 -\frac{t}{z}J_1 + tzJ_1},
\quad z\to 0.
\end{eqnarray}
Now the jump matrix for $\Psi(z)$ is a constant matrix :
\begin{equation}\label{e-vconst}
\Psi_+(z)=\Psi_-(z)
 \begin{pmatrix} 0&1\\-1&1 \end{pmatrix}  , \qquad z\in\Sigma.
\end{equation}

{\bf Equations :}
By differentiating \eqref{e-vconst},
$\frac{\partial}{\partial z} \Psi$ satisfies the same jump condition.
Hence $\frac{\partial\Psi}{\partial z}  \Psi^{-1}$ has no jump cross
$\Sigma$.
Also from \eqref{e-Psi}, $\frac{\partial\Psi}{\partial z}  \Psi^{-1}$
has a double pole at $z=0$.
Therefore using \eqref{e-Psias1}, \eqref{e-Psias2}, we have
\begin{equation}\label{e-eqz}
  \frac{\partial\Psi(z;k)}{\partial z}  \Psi(z;k)^{-1}
= tJ_1 + \frac{A_1(k)}{z} + \frac{A_2(k)}{z^2}=:P(z;k),
\end{equation}
for some constant matrices $A_1, A_2$ which depend on $k$.
On the other hand, $\Psi(z;k+1)\Psi(z;k)^{-1}$ again has no
jump cross $\Sigma$, and it is now entire.
Thus using \eqref{e-Psias1}, \eqref{e-Psias2}, we have
\begin{equation}\label{e-eqk}
  \Psi(z;k+1)\Psi(z;k)^{-1}= zJ_1 + X_0(k) =:Q(z;k),
\end{equation}
for some constant matrix $X_0(k)$.

Now we take `cross differentiation' of \eqref{e-eqz} and \eqref{e-eqk}
in two different ways :
\begin{equation}
\begin{split}
  \frac{\partial}{\partial z} \Psi(z;k+1)
&=  \frac{\partial}{\partial z} Q(z;k)\Psi(z;k)
= \frac{\partial Q(z;k)}{\partial z}\Psi(z;k)
+ Q(z;k)\frac{\partial \Psi(z;k)}{\partial z} \\
&= \biggl\{ \frac{\partial Q(z;k)}{\partial z} + Q(z;k)P(z;k) \biggr\}
\Psi(z;k) .
\end{split}
\end{equation}
and
\begin{equation}
\begin{split}
\frac{\partial \Psi(z;k)}{\partial z}
 = \frac{\partial \Psi(z;k)}{\partial z}\biggr|_{k\mapsto k+1}
&= P(z;k+1) \Psi(z;k+1)
= P(z;k+1)Q(z;k) \Psi(z;k).
\end{split}
\end{equation}
Thus we obtain an equation
\begin{equation}
\frac{\partial Q(z;k)}{\partial z} + Q(z;k)P(z;k)
=  P(z;k+1)Q(z;k) .
\end{equation}
By plugging in the formulas of $P$ and $Q$ and comparing the coefficients
in $z$, we obtain the relations for $A_1, A_2, X_0$ :
\begin{eqnarray}
\label{e-e1}  &&J_1+J_1A_1(k)-A_1(k+1)J_1 + t[X_0, J_1]=0, \\
\label{e-e2}  &&J_1A_2(k)-A_2(k+1)+X_0(k)A_1(k)-A_1(k+1)X_0(x)=0, \\
\label{e-e3}  &&X_0(k) A_2(k)-A_2(k+1)X_0(k)= 0.
\end{eqnarray}

{\bf Constant matrices $A_1, A_2, X_0$ : }
Now we express $A_1, A_2, X_0$ in terms of $m$.
We plug in \eqref{e-Psias1} and \eqref{e-Psias2} into
\eqref{e-eqz} and \eqref{e-eqk}.
This determines the constant matrices $A_1, A_2, X_0$.
Using \eqref{e-Psias1} for \eqref{e-eqz},
\begin{eqnarray}
\label{e-A1-1} O(z^{-1}) :  A_1(k) = t[m^\infty_1(k), J_1] + kJ_1.
\end{eqnarray}
Using \eqref{e-Psias2} for \eqref{e-eqz},
\begin{eqnarray}
\label{e-A1-2} &O(z^{-1})& :  A_1(k) = B(k)\bigl\{ t[m^0_1(k), J_1] + kJ_1
\bigr\} B(k)^{-1}  , \\
\label{e-A2} &O(z^{-2})& : A_2(k)= tB(k) J_1 B(k)^{-1}.
\end{eqnarray}
Using \eqref{e-Psias1} for \eqref{e-eqk},
\begin{eqnarray}
\label{e-X-1} O(1) :  X_0(k) = J_0 + m^\infty_1(k+1)J_1-J_1 m^\infty_1(k).
\end{eqnarray}
Using \eqref{e-Psias2} for \eqref{e-eqk},
\begin{eqnarray}
\label{e-J} &O(z)& : J_1=  B(k+1)\bigl\{ J_1 + m^0_1(k+1)J_0-J_0m^0_1(k)
\bigr\} B(k)^{-1}  ,\\
\label{e-X-2} &O(1)& :  X_0(k) = B(k+1) J_0 B(k)^{-1}.
\end{eqnarray}

{\bf Symmetry :}
Note that the jump matrix $v(z)$ for $m$ satisfies
\begin{equation}
  \sigma_3 v(z^{-1})^t \sigma_3= v(z),
\qquad \sigma_3:= \begin{pmatrix} 1&0\\0&-1 \end{pmatrix}.
\end{equation}
Thus it is standard to show that $m$ satisfies the symmetry condition
\begin{equation}
  m(z)= \sigma_3 B^t(m(z^{-1})^t)^{-1}\sigma_3.
\end{equation}
By taking $z\to 0$, and considering terms of order $O(1)$ and $O(z)$,
we obtain the symmetry relations
\begin{eqnarray}
   B&=& \sigma_3B^t\sigma_3, \\
  m^0_1&=& -\sigma_3 (m^\infty_1)^t\sigma_3.
\end{eqnarray}
Since $\det B=1$ and $tr(m^0_1)=0$, we can set
\begin{equation}\label{e-B}
  B(k):=\begin{pmatrix} a(k)&b(k)\\-b(k)&d(k) \end{pmatrix} ,
\qquad a(k)d(k)+b(k)^2=1,
\end{equation}
and set
\begin{equation}\label{e-m}
  m^0_1(k)
:=\begin{pmatrix} \lambda(k)&\mu(k)\\ \nu(k)&-\lambda(k) \end{pmatrix} ,
\qquad   m^\infty_1(k)
:=\begin{pmatrix} -\lambda(k)&\nu(k)\\  \mu(k)&\lambda(k) \end{pmatrix}.
\end{equation}

\bigskip

{\bf Equations from \eqref{e-e1}-\eqref{e-e3} : }
Now we re-write the equations \eqref{e-e1}-\eqref{e-e3} in terms
$a,b,d$ and $\lambda, \mu, \nu$.
First, if we use \eqref{e-A1-1} and \eqref{e-X-2},
use \eqref{e-B} and \eqref{e-m}, the equation \eqref{e-e1} yields
\begin{eqnarray}
\label{e-mbd} \mu(k)&=&b(k+1)d(k), \\
\label{e-nab} \nu(k+1)&=& -a(k+1)b(k).
\end{eqnarray}
Second, if we use \eqref{e-A1-2}, \eqref{e-A2} and \eqref{e-B}, \eqref{e-m},
the equation \eqref{e-e2} yields the same equations \eqref{e-mbd} and
\eqref{e-nab}.
Finally, the equation
\eqref{e-e3} is trivial if we use \eqref{e-A2} and \eqref{e-X-2}.

{\bf Equations from \eqref{e-A1-1}-\eqref{e-J} : }
We also have the condition that the two different formulas
\eqref{e-A1-1} and \eqref{e-A1-2} for $A_1(k)$ are equal.
Similarly \eqref{e-X-1} and \eqref{e-X-2} for $X_0(k)$
are equal. Also there is an additional condition \eqref{e-J}.
From \eqref{e-A1-1} and \eqref{e-A1-2}, we have
\begin{equation}
  \bigl( t[m^\infty_1(k), J_1] + kJ_1 \bigr) B(k)
= B(k) \bigl( t[m^\infty_1(k), J_1] + kJ_1 \bigr).
\end{equation}
Using \eqref{e-B}, \eqref{e-m}, this yields
\begin{equation}\label{e-all}
  t\bigl(\mu(k)a(k)-\nu(k)d(k)\bigr)+kb(k)=0.
\end{equation}
If we use equations \eqref{e-mbd} and \eqref{e-nab}, and the relation
$ad+b^2=1$, this equation becomes
\begin{equation}
  kb(k)+t\bigl(b(k-1)+b(k+1)\bigr)(1-b(k)^2))=0.
\end{equation}

From \eqref{e-X-1} and \eqref{e-X-2}, we have
\begin{equation}
  \bigl( J_0+ m^\infty_1(k+1)J_1- J_1m^\infty_1(k) \bigr) B(k)
= B(k+1) J_0,
\end{equation}
which yields
\begin{eqnarray}
\label{e-a1}  a(k+1)-a(k)+\nu(k+1)b(k) &=& 0, \\
 b(k)+\nu(k+1) d(k) &=& 0, \\
b(k+1) -(\lambda(k+1)-\lambda(k))b(k)-\mu(k)a(k) &=&0,\\
(\lambda(k+1)-\lambda(k))d(k)-\mu(k)b(k) &=&0.
\end{eqnarray}
Thus using \eqref{e-nab}, \eqref{e-a1} yields
\begin{equation}
  a(k)= \bigl(1-b(k)^2\bigr) a(k+1).
\end{equation}

On the other hand, the equation \eqref{e-J} implies
\begin{equation}
  J_1B(k)
= B(k+1) \bigl( J_1+ m^0_1(k+1)J_0- J_0m^0_1(k) \bigr),
\end{equation}
which yields
\begin{eqnarray}
(\lambda(k+1)-\lambda(k))a(k+1)-\nu(k+1)b(k+1) &=&0, \\
 b(k+1)-\mu(k) a(k+1) &=& 0, \\
b(k) -(\lambda(k+1)-\lambda(k))b(k+1)-\nu(k+1)d(k+1) &=&0,\\
\label{e-d1} d(k+1)-a(k)+\mu(k+1)b(k) &=& 0.
\end{eqnarray}
Thus using \eqref{e-mbd}, \eqref{e-d1} yields
\begin{equation}
  d(k)= \bigl(1-b(k)^2\bigr) d(k-1).
\end{equation}

\end{proof}

\begin{rem}
The RHP \eqref{e-rhp} has also the parameter $t$. In addition to the above
equation \eqref{e-diffeq}-\eqref{e-diffeqd} obtained
from the $z$-derivative and
`$n$-derivative' of $Y$, we can also obtain two more equations from
$z$-derivative and $t$-derivative of $Y$, and from the $t$-derivative and
$n$-derivative of $Y$.
In this paper, we would not consider those equations.
See Section 3 of \cite{BDJ} and \cite{TW3, Wang, AvM}
for those other equations.
\end{rem}

\section{Integrable operators}\label{sec-intg}

In this section, we related the RHP \eqref{e-rhp} to a so-called integrable
operator on $L^2(\Sigma)$, and also discuss an identity between
Toeplitz determinants and Fredholm determinants on $L^2(\Sigma)$.

It is shown in \cite{BDJ2} that
the RHP \eqref{e-rhp} for
orthogonal polynomials is also the RHP for an integrable operator on $\Sigma$
after a simple algebraic transformation.
For this purpose, the particular
choice $\varphi(z)=e^{t(z+z^{-1})}$ in \eqref{e-rhp} is not important.
In the below, we just assume that $\varphi(z)$ is a continuous
function on $\Sigma$ which has the factorization
$\varphi(z)=\varphi_+(z)\varphi_-(z)$ where $\varphi_+(z), \varphi_+(z)^{-1}$
are analytic in $|z|<1$, continuous in $|z|\le 1$, $\varphi_+(0)=1$,
and $\varphi_-(z), \varphi_-(z)^{-1}$ are analytic in $|z|>1$,
continuous in $|z|\ge 1$, $\lim_{z\to\infty}\varphi_-(z)=1$.
We also assume that $1$ is not in the spectrum
of the operator $K$ defined in \eqref{e-intop} below.

Define (cf. \eqref{e-mdef1}, \eqref{e-mdef2})
\begin{equation}\label{e-mdef11}
  M(z):= \begin{pmatrix} 0&-1\\ 1&0 \end{pmatrix} Y(z)
\begin{pmatrix} \frac12 \varphi_+(z)&0\\0& 2\varphi_+(z)^{-1} \end{pmatrix}  ,
\qquad |z|<1,
\end{equation}
and set
\begin{equation}\label{e-mdef22}
  M(z):= \begin{pmatrix} 0&-1\\1&0 \end{pmatrix} Y(z)
\begin{pmatrix} 0&z^{-k}\varphi_-(z)\\ -z^k\varphi_-(z)^{-1}&0\end{pmatrix},
\qquad |z|>1.
\end{equation}
Then $M(z)$ solves the new Riemann-Hilbert problem
\begin{equation}\label{e-newrhp2}
\begin{cases}
  M(z) \qquad \text{is analytic in $z\in \Omega_\pm$,
and is continous in $\overline{\Omega_\pm}$,}\\
  M_+(z)=M_-(z) V(z), \qquad z\in\Sigma,\\
  M(z) = I +O(z^{-1}), \qquad z\to\infty,
\end{cases}
\end{equation}
where
\begin{equation}
V(z)= \begin{pmatrix} 0& - 2z^{-k}\psi(z)^{-1}\\
\frac12 z^{k}\psi(z)&2  \end{pmatrix},
\qquad \psi(z):= \varphi_+(z)\varphi_-(z)^{-1}.
\end{equation}
Note that the jump matrix $V$ is of the form
\begin{equation}\label{e-intjump}
 V(z)=I-2\pi i f(z)g(z)^T,
\end{equation}
where
\begin{equation}
 f(z)=(f_1, f_2)^T = (z^{-k}, -\frac12\psi(z))^T,
\qquad g(z)=(g_1, g_2)^T= \frac1{2\pi i} (z^k, 2\psi(s)^{-1})^T.
\end{equation}
This is precisely the RHP
canonically associated to the
so-called integrable operator $K=K_k$ on $L^2(\Sigma)$
(see \cite{IIKS}, \cite{deiftint})
whose kernel is defined by
\begin{equation}\label{e-intop}
  K(z,w)= \frac{f(z)^Tg(w)}{z-w}
= \frac{z^{-k}w^k-\psi(z)\psi(w)^{-1}}{2\pi i(z-w)},
\qquad (Kf)(z)=\int_{\Sigma} K(z,w)f(w)dw.
\end{equation}

It is shown in \cite{IIKS} (see also \cite{deiftint})
that given an integrable operator $K$ of the form \eqref{e-intop}
(with general column vectors $f$, $g$ such that $f(z)^Tg(z)=0$),
the resolvent kernel of $K$, if exists, is again an integrable
operator
\begin{equation}
  \biggl( \frac{1}{1-K}K \biggr) (z,w) = \frac{F(z)^TG(z)}{z-w}
\end{equation}
and moreover,
the column vector $F$, $G$ are expressed in terms of the solution $M$
to the (normalized) RHP with the jump matrix \eqref{e-intjump} :
\begin{eqnarray}
  F(z) = \biggl( \frac1{1-K}f \biggr)(z) &=& M_+(z)f(z), \\
  G(z) = \biggl( \frac1{1-K}g \biggr)(z) &=& (M_+(z)^T)^{-1}g(z).
\end{eqnarray}

From the structure of the
integrable operators, one can also show (see Lemma 4 of \cite{BDJ2}) that
\begin{equation}\label{e-M11}
  M_{11}(0) = \frac{\det(1-K_{k-1})}{\det(1-K_k)}.
\end{equation}
Thus, using the definition of $M$, \eqref{e-N} and \eqref{e-ND},
we obtain
\begin{equation}\label{e-detratio}
  \frac{D_{k-1}(\varphi)}{D_k(\varphi)}
= 2 \frac{\det(1-K_{k-1})}{\det(1-K_k)}.
\end{equation}
Moreover, when $[1,\infty)$ has no intersection to the spectrum of $K_k$
for all $k\ge k_0$ for some $k_0$, it is shown in Lemma 5 of \cite{BDJ2}
that
\begin{equation}
  \lim_{p\to\infty} 2^{-p} \det(1-K_p)= 1.
\end{equation}
Thus for this case, by taking infinite product in $k$ of \eqref{e-detratio},
we have
\begin{equation}\label{e-detinf}
  \frac{D_n(\varphi)}{D_\infty(\varphi)} = 2^{-n}\det(1-K_n),
\end{equation}
where
\begin{equation}
  D_\infty(\varphi)=\lim_{p\to\infty} D_p(\varphi)
= \exp(\sum_{j=1}^\infty j(\log \varphi)_j (\log \varphi)_{-j})
\end{equation}
is given by the strong Szeg\"o theorem.

It would be of interest to find if
the other entries $M_{12}$, $M_{21}$ have expressions involving $K$,
analogous to \eqref{e-M11}.
Such formulas, if exist, would yield similar formulas to
\eqref{e-detratio}, \eqref{e-detinf} for the Hankel determinants
\eqref{e-Haneklell}.

We note that as mentioned in Remark \ref{rem-daku},
there is yet another identity between the Toeplitz determinant
and the Fredholm determinant of an integrable operator.
But in this case, the integrable operator acts on a discrete set
$\{k,k+1,\cdots\}$.
Hence there is another RHP with jump conditions on the discrete set.
It is not clear yet
if for example, there is an algebraic transformations between
these two RHP's, one with jump on $\Sigma$ and the other
with jump on the discrete set.

\section{Other models}\label{sec-othe}

There are other last passage percolation models which have
the similar Toeplitz/Hankel determinant formulas of the form \eqref{e-U},
\eqref{e-O} and \eqref{e-extform}.
The only difference now is the symbol $\varphi$ which depends
on the model.
In the below, we will describe the models and the corresponding symbols.
Here we only consider the square case and the triangle case.

Since the RHP formulation for the Toeplitz/Hankel determinants
in Section \ref{sec-perc} is general,
for all the models below there are associated RHP
of the form \eqref{e-rhp} with new function $\varphi(z)$.
Thus we can follows the procedures in Section \ref{sec-asym},
Section \ref{sec-diff} and Section \ref{sec-intg}.
The determinant identity
in Section \ref{sec-intg} is general and it includes all
the cases below.
The difference equations for the cases below can be worked out as in Section
\ref{sec-diff}.
It will appear in some other place.
See \cite{AvM}
for some results for the difference/differential equations.
The asymptotic analysis of Section \ref{sec-asym} can be worked out
in all cases below.
However, as mentioned several places above,
for the square case,
the convergence in distribution result of the form \eqref{e-FGUE}
can be obtained from a Fredholm determinant formula on a discrete set,
which does not require the RHP analysis.
In the literature, many such results are obtained from the Fredholm
determinant analysis.
But for the triangle cases and also for the tail estimates,
RHP analysis have been used (see \cite{B00, BR3}).
For the survey of asymptotic results of the cases below,
see e.g. \cite{Rains:mean}.

\subsection{Square case}

Consider a lattice version of the Poisson process.
Let $X(i,j)$, $i,j\in \N$, be
a planar array of independent random variables.
We consider a \emph{weakly-up/weakly-right
path}, which a sequence $\{(i_k,j_k)\}_{k=1}^l$
such that $i_1\le i_2\le\cdots\le i_l$ and
$j_1\le j_2\le \cdots\le j_l$.
Let $(1,1)\nearrow(M,N)$ denote the set of all
weakly-up/weakly-right paths from $(1,1)$ to $(M,N)$.
Now define
\begin{equation}\label{e-i3}
 L(M,N):= \max_{\pi\in (1,1)\nearrow(M,N)} \sum_{(i,j)\in\pi} X(i,j).
\end{equation}
If one thinks of the random variable
$X(i,j)$ as the passage time at the site $(i,j)$, the
sum $\sum_{(i,j)\in\pi} X(i,j)$ is
the total passage time to travel from $(1,1)$ to $(M,N)$
along the particular path $\pi$.
Thus $L(M,N)$ can be regarded as the \emph{last passage
percolation time} from $(1,1)$ to $(M,N)$.
As in the Poisson case of Section \ref{sec-perc},
this problem arises in many different fields.
See for example, \cite{kurtj:shape, GTW, SpohnP3}
for various applications.

In the definition of $L(M,N)$ in \eqref{e-i3}, we could take different admissible up/right
paths. For example, we could take weakly-up/strictly-right paths (sequences
$\{(i_k,j_k)\}_{k=1}^l$ such that $i_1\le\cdots\le i_l$, $j_1<\cdots<j_l$) or
strictly-up/strictly-right paths (sequence $\{(i_k,j_k)\}_{k=1}^l$ such that $i_1<\cdots<
i_l$, $j_1<\cdots<j_l$). Also we can take continuum model instead of lattice model. In the
below, we list the models which have the determinant formula for the distribution of
$L(M,N)$. In each case, we have
\begin{equation}\label{e-Uform}
  \Prob(L(M,N)\le \ell)
= \frac{1}{Z_{M,N}} \Exp_{U\in U(\ell)} \det(\varphi(U)).
\end{equation}
with a function $\varphi(z)$ for $|z|=1$ which depends on the model.
The constant $Z_{M,N}=\lim_{\ell\to\infty} D_\ell(\varphi)$ is
a finite number in each case.

We use the notation $g(q)$ for a geometric random variable with
parameter $q\in (0,1)$ :
$\Prob(g(q)=k)=(1-q)q^k$, $k=0,1,2,\cdots$.
We understand that $g(0)$ means the identically $0$ random variable.
The notation $b(q)$
is used for the Bernoulli
random variable : $\Prob(b(q)=0)=\frac{q}{1+q}$, $\Prob(b(q)=1)=\frac1{1+q}$.

{\bf Lattice models} :
\begin{itemize}
\item[(a)] (weakly-up/weakly-right) Fix $q_i, q_j\ge 0$ such that
$q_iq_j \in [0,1)$, $i=1,\cdots,M$, $j=1,\cdots, N$.
Take $X(i,j)\sim g(q_iq_j)$. Then
\begin{equation}\label{}
  \varphi(z)= \prod_{i=1}^M\prod_{j=1}^N (1+q_iz)(1+q_j/z),
  \qquad Z_{M,N}=\prod_{i=1}^M\prod_{j=1}^N (1-q_iq_j)^{-1}.
\end{equation}

\item[(b)] (weakly-up/strictly-right) Fix $q_i, q_j\ge 0$,
$i=1,\cdots,M$, $j=1,\cdots, N$.
Take $X(i,j)\sim b(q_iq_j)$. Then
\begin{equation}\label{}
  \varphi(z)= \prod_{i=1}^M\prod_{j=1}^N (1+q_iz)(1-q_j/z)^{-1},
  \qquad Z_{M,N}=\prod_{i=1}^M\prod_{j=1}^N(1+q_iq_j).
\end{equation}

\item[(c)] (strictly-up/strictly-right) Fix $q_i, q_j\ge 0$
such that $q_iq_j\in[0,1)$, $i=1,\cdots,M$, $j=1,\cdots, N$.
Take $X(i,j)\sim g(q_iq_j)$. Then
\begin{equation}\label{}
  \varphi(z)= \prod_{i=1}^M\prod_{j=1}^N (1-q_iz)^{-1}(1-q_j/z)^{-1},
  \qquad Z_{M,N}=\prod_{i=1}^M\prod_{j=1}^N(1-q_iq_j)^{-1}.
\end{equation}

\end{itemize}

{\bf Lattice-Continuum models} :

\begin{itemize}

\item[(d)] (Poisson : weakly-right) Fix $q_i\ge 0$, $i=1,\cdots, N$.
Consider a Poisson process of rate $q_i$ on
$\R\times \{i\}\subset \R\times \{1,\cdots,N\}$.
Let $L(t,N)$ be the length of the longest weakly-up/weakly-right path
from $(0,1)$ to $(t,N)$. Then
$\Prob(L(t,N)\le \ell)=\frac1{Z_{t,N}} \Exp_{U\in U(\ell)}(\varphi)$ with
\begin{equation}\label{}
  \varphi(z)= e^{tz}\prod_{i=1}^N (1+q_i/z),
  \qquad Z_{t,N}=\prod_{i=1}^N e^{t q_i}.
\end{equation}

\item[(e)] (Poisson : strictly-right)
In (f), replace the admissible path as weakly-up/strictly-right paths. Then
\begin{equation}\label{}
  \varphi(z)= e^{tz}\prod_{i=1}^N (1-q_i/z)^{-1},
  \qquad Z_{t,N}=\prod_{i=1}^N e^{t q_i}.
\end{equation}

\end{itemize}

{\bf Continuum model} :

\begin{itemize}

\item[(f)] (Poisson process)
This is the model in the square case of Section \ref{sec-perc}.
We have
$\Prob(L(t)\le \ell)=\frac1{Z_{t}} \Exp_{U\in U(\ell)}(\varphi)$ with
\begin{equation}\label{}
  \varphi(z)= e^{t(z+z^{-1})},
  \qquad Z_{t}=e^{-t^2}.
\end{equation}

\end{itemize}

\begin{rem}
In (d) and (e), weakly-up condition can be changed to strictly-up, since in
a Poisson process, the event of two points occurring at the same position
has measure $0$.
\end{rem}

\begin{rem}
We could have mixture of the above models. Then the corresponding $\varphi$
is simply multiplication of each symbol. See \cite{Rains:mean} for the most
general setting.
\end{rem}

\subsection{Symmetrized models}

For models (a), (c), (f), one can think of 4 different symmetrized models
(see section 7 of \cite{BR1}). We could impose that the model is symmetric
under (1) rotation about the center, (2) reflection about the diagonal,
(3) reflection about the anti-diagonal and (4) reflection about both
the diagonal and the anti-diagonal.
In each case, there are still determinant formulas for the distribution of
the longest up/right path.
However, for the symmetry types (2),(3),(4) for (a), (c),
we only have formulas for the square model $M=N$.
Depending on the symmetry type,
instead of Toeplitz determinant above, we sometimes have Hankel determinant.

Here we will consider only the case (2).
That is, for example in (a), we impose that $X(j,i)$ should be equal to the
value of $X(i,j)$ for each $i,j$.
In each case, the algebraic formula is of the form
\begin{equation}\label{e-Oform}
  \Prob(L\le \ell) = \frac1{Z} \Exp_{U\in O(\ell)} \det(\psi(U)).
\end{equation}
Note that now we have expectation over the orthogonal group.
We will specify the function $\psi(z)$.

Let $g'(\alpha, q)$ be the random variable given by
\begin{equation}
  \Prob(g'(\alpha, q)=k)= \frac{1-q^2}{1+\alpha q} \alpha^{\text{k mod 2}}q^k.
\end{equation}

\begin{itemize}
\item[(a-S)] Fix $q_i\ge 0$, $i=1,\cdots, N$ and $\alpha\ge 0$
such that $q_i\in[0,1)$ and $\alpha q_i\in[0,1)$.
Take $X(i,j)=X(j,i)\sim g(q_iq_j)$ for $i\neq j$ and
$X(i,i) \sim g(\alpha q_i)$. Then the longest weakly-up/weakly-right
$L(N;\alpha)$
path from $(1,1)$ to $(N,N)$ satisfies \eqref{e-Oform} with
\begin{equation}
  \psi(z)= (1+\alpha z) \prod_{i=1}^N (1+q_i z),
\qquad Z= \prod_{i=1}^N (1-\alpha q_i)^{-1}
\prod_{1\le i<j\le N} (1-q_iq_j)^{-1}.
\end{equation}

\item[(c-S)]
Fix $q_i\ge 0$, $i=1,\cdots, N$ and $\alpha\ge 0$
such that $q_i\in[0,1)$.
Take $X(i,j)=X(j,i)\sim g(q_iq_j)$ for $i\neq j$ and
$X(i,i) \sim g'(\alpha , q_i)$. Then the longest strictly-up/strictly-right
$L(N;\alpha)$
path from $(1,1)$ to $(N,N)$ satisfies \eqref{e-Oform} with
\begin{equation}
  \psi(z)= (1+\alpha z) \prod_{i=1}^N (1-q_i z)^{-1},
\qquad Z= \prod_{i=1}^N (1+\alpha q_i)(1-q_i^2)^{-1}
\prod_{1\le i<j\le N} (1-q_iq_j)^{-1}.
\end{equation}

\item[(f-S)] The Triangle case of Section \ref{sec-perc}.
Take a Poisson process of rate $1$ in the region $x<y$ in the
xy-plane. Then take the mirror image of the points about the line $x=y$.
In addition, take a Poisson process of rate $\alpha\ge 0$ on the
line $x=y$. Then the longest up/right path $L(t;\alpha)$ from $(0,0)$
to $(t,t)$ satisfies \eqref{e-Oform} with
\begin{equation}
  \psi(z)= (1+\alpha z) e^{tz}, \qquad Z= e^{\alpha t + t^2/2}.
\end{equation}

\end{itemize}



\bibliographystyle{plain}
\bibliography{main}

\end{document}

%% file: newcommands.tex
\DeclareMathOperator{\Exp}{\mathbb E}

\DeclareMathOperator{\Prob}{\mathbb P}

\newcommand{\C}{\mathbb C}
\newcommand{\R}{\mathbb R}

\newcommand{\N}{\mathbb N}